\documentclass[12pt,a4paper]{article}

\usepackage[latin1]{inputenc}
\usepackage{amsmath, amsfonts, amssymb, graphicx, amsthm}
\usepackage{graphicx}
\usepackage[spaces,hyphens]{url}
\makeindex
\newtheorem{theorem}{Theorem}[section]
\newtheorem{corollary}[theorem]{Corollary}
\newtheorem{defi}[theorem]{Definition}

\newtheorem{lemma}[theorem]{Lemma}

\newtheorem{proposition}[theorem]{Proposition}
\newtheorem{remark}[theorem]{Remark}

\newcommand*{\N}{\mathbb{N}}

\newcommand*{\C}{\mathbb{C}}
\newcommand*{\Bcal}{\mathcal{B}}

\newcommand*{\tr}{\text{{\normalfont{tr}}}\,}

\newcommand*{\1}{\text{{\sf 1}}}
\newcommand{\bra}[1]{\langle\,#1\,\rvert}
\newcommand{\ket}[1]{\lvert\,#1\,\rangle}

\newcommand{\norm}[1]{\left\Arrowvert#1\right\Arrowvert}
\newcommand{\abs}[1]{\left\arrowvert#1\right\arrowvert}

\newcommand{\inner}[2]{\langle\,#1\,,\,#2\,\rangle}

\newcommand{\dm}{\,d\mu}
\newcommand{\dn}{d\nu}
\renewcommand{\Pr}{\mathbb{P}}
\newcommand*{\E}{\mathbb{E}}
\newcommand*{\supp}{\text{supp}\,}

\newcommand{\rank}{\text{rank}\,}
\newcommand{\npurifies}{\,\,n\text{-purifies}\,\,}
\newcommand{\biggnorm}[1]{\bigg\Arrowvert #1 \bigg\Arrowvert}
\newcommand{\im}{Im\,\,}

\title{ Lyapunov exponents for Quantum Channels: an entropy formula and generic properties }

\author{Jader E. Brasil, Josu\'e Knorst and Artur O. Lopes}

\begin{document}
\maketitle

\begin{abstract}
We denote by  $M_k$ the set  of  $k$ by $k$ matrices with complex entries.
We consider quantum channels $\phi_L$ of the form: given a measurable function   $L:M_k\to M_k$ and a measure   $\mu$  on $M_k$  we define the linear operator $\phi_L:M_k \to M_k$, by the law $\rho \,\to\,\phi_L(\rho) = \int_{M_k} L(v) \rho L(v)^\dagger \, \dm(v).$

In a previous work, the authors show that for a fixed measure $\mu$ the $\Phi$-Erg property is generic on the function $L$  (also irreducibility). Here we will show that the purification property is also generic on $L$ for a fixed $\mu$.

Given $L$ and $\mu$ there are two related stochastic processes:  one takes values on the projective space $ P(\C^k)$ and the other  on matrices in $M_k$.
The   $\Phi$-Erg property and the purification condition are the nice hypothesis for the discrete time evolution given by the natural transition probability. In this way it will follow that generically on  $L$,
if  $\int |L(v)|^2 \log |L(v)|\, \ d\mu(v)<\infty$,  the Lyapunov exponents $\infty > \gamma_1\geq \gamma_2\geq ...\geq \gamma_k\geq -\infty$ are well defined.

In a previous work, the concepts of entropy of a channel and Gibbs channel were presented; and  also an example (associated to a stationary  Markov chain) in which this definition of entropy (for a quantum channel)  matches the Kolmogorov-Shanon definition of entropy.  We estimate here the larger Lyapunov exponent for the mentioned example and we show that it is equal to $-\frac{1}{2} \,h$, where $h$ is the entropy of the associated  Markov invariant probability.
\end{abstract}

 \vspace{1cm} \maketitle

%\begin{document}
	\section{Introduction}

We consider quantum channels of the form $\phi_L:M_k \to M_k$,  in which $M_k$ is the set  of  complex $k$ by $k$ matrices, $\rho \,\to\,\phi_L(\rho) = \int_{M_k} L(v) \rho L(v)^\dagger \, \dm(v)$, $\,L:M_k\to M_k$ is a measurable function and $\mu$ is a measure on $M_k$.

In a recent work \cite{benoist2017invariant},  the authors considered Lyapunov exponents for this class of channels  $\phi_L$,  when $L$ was constant and equal to the identity matrix. The $\Phi$-Erg property and the purification condition (see definitions on section \ref{gen})  were also considered in this mentioned paper.

In a previous paper (see \cite{BKL}), we show  that for a fixed measure $\mu$ it is generic on the function $L$ the $\Phi$-Erg property (in fact we show that the irreducible condition is generic). The novelty here is that we will show that the purification condition is also generic on $L$ for a fixed measure $\mu$ (see section \ref{purgen}).

The introduction of this variable $L$ allows us to consider questions of  generic nature in this type of problem. We use the $C^0$ topology in the set of complex matrices.

For the benefit of the reader in the Appendix  Section \ref{appe} we present an outline of the results in \cite{BKL} and  relations of Lyapunov exponents with previous works. 

Following  \cite{benoist2017invariant}, one can consider
associated to $L$ and $\mu$  two related processes: one denoted by  $X_n$, $n\in \mathbb{N}$, takes values on the projective space $ P(\C^k)$; and the other,  denoted by $\rho_n$, $n\in \mathbb{N}$, takes values on $D_k$  (where $D_k$ is
the set of density operators). The natural transition probability is defined in \cite{benoist2017invariant}.

The   $\Phi$-Erg property and the purification property play an important role when analyzing the ergodic properties of these two processes (see Section \ref{gen}).

 Here we consider a concept of quantum entropy for channels in section \ref {exam} which was 
 initially presented in \cite{baraviera2010thermodynamic}.
 
For a fixed $\mu$ and a general $L$, it was  presented in \cite{BKL}  a natural concept of entropy  (see   future section \ref{ent})  for a channel in order to develop a version of  Gibbs formalism in this setting. It was also presented in example 8.5 in \cite{BKL} a certain channel (related to stationary Markov Chains)  in which  the value obtained with  this definition coincides with the classical  value of entropy.  This shows that the concept that was introduced is natural. This definition  of entropy is a generalization of the concept described on papers \cite{baraviera2010thermodynamic}, \cite{BLLT2} and \cite{BLLT1}.
This particular form of defining entropy is, in a certain way, inspired by the results of
\cite{slomczynski2003dynamical}  which considers iterated function systems.

We call the example described in example 8.5 in \cite{BKL} 
the {\it Markov model in quantum information}. It is the main example to be considered in our Section \ref{exam}.

 The main contribution to the topic of Lyapunov exponents of quantum channels   is the paper \cite{benoist2017invariant}.
 We adapt here  the formalism of \cite{benoist2017invariant}  to the case of a general $L$
in order to estimate the Lyapunov exponents of the associated  dynamical time evolution. In this part of the paper, we  just outline the main points of the  proof (see section 7)
because it is basically the one presented  in \cite{benoist2017invariant}.
We describe the sufficient conditions for the  Lyapunov exponents to be finite. The irreducibility and the purification condition are the appropriate assumptions (as mentioned in \cite{benoist2017invariant}). Proposition \ref{phi_erg_dense_set1} and \cite{BKL} show that those conditions are true for a fixed $\mu$ and a generic $L$.

Results relating entropy and Lyapunov exponents (both in the classical sense) are quite important in Ergodic Theory (see the Appendix section  \ref{appe} or \cite{Bar}, \cite{ManS}, \cite{Mane} and \cite{Kat}).

We compute the first Lyapunov exponent (which is negative) for the  {\it Markov model in quantum information} (see Section \ref{exam}) and we show that it is equal to $-\frac{1}{2} \,h$, in which $h$ is the entropy of the associated  Markov probability. We also show that the second Lyapunov exponent, in this case,  will be
$-\infty$. Naturally, a similar general result for the  class of all quantum channels is not reachable due to its inherent  generality.

We point out that the definition of entropy for a (normalized) channel presented in \cite{BKL} explores the use of a ``kind'' of Ruelle operator. This  procedure uses a natural {\it a priori} probability   and this is plausible  due to the fact that the ``set of preimages'' can be uncountable (see \cite{LMMS}). The main issue on the reasoning in \cite{BKL} is invariance (in time one); however, in that paper, the concept of entropy is not directly related (we mean: in an explicit form) to time evolution.
On the other hand, we point out that the values of the Lyapunov exponents are of {\it dynamical nature}.  The dynamical discrete-time evolution - on the present setting - is described by a stochastic process taking values on the set of matrices in $M_k$ (see sections \ref{gen} and \ref{Ly}). The example we consider in section \ref{exam} shows that  the concept of entropy of a channel (at least in this case) presented in \cite{BKL} can be linked to the natural dynamical time evolution via the main Lyapunov exponent. In the Appendix Section \ref{appe} we present the relation of our results about the Lyapunov exponents with the analogous  classical concept in Dynamical Systems.

General references for Quantum Channels and Quantum Mechanics can be found at \cite{Li}, \cite{Petz2}, \cite{Alicki}, \cite{baraviera2010thermodynamic}, \cite{CFF},  \cite{BLMM}, \cite{LCol}, \cite{LSeb}, \cite{LSeb0} and \cite{wolf2012quantum}. The book
\cite{Ma1} presents several important results for the general theory of Lyapunov exponents (see \cite{TV}, \cite{DK1}, \cite{Bac},  and \cite{DK}). \cite{Ma} and \cite{PP} describe basic result in Ergodic Theory.

We thank S. Klein for supplying us with references and also Cyntia Trevor who helped us in the grammatical writing of the text.

\section{Basic results} \label{bas}

In this section,  for the benefit of the reader, we recall some well-known preliminary results  on the topic.
We will mention throughout the section the references where the different results are demonstrated.

We denote by $M_k$ the set of complex $k$ by $k$ matrices.
We denote by $ \text{Id}$ the identity matrix on $M_k$.

 We consider the standard Borel sigma-algebra over $M_k$ and the canonical Euclidean inner product on $\mathbb{C}^k$

According to our notation, $\dagger$ denotes the operation  of taking the dual of a matrix with respect to the canonical inner product on $\C^k$.

			Here tr denotes the trace of a matrix.

	Given two matrices $A$ and $B$ we define the Hilbert-Schmidt product
	$$\inner{A}{B} \,\,= \, \text{ tr}\,\,  (A\, B^\dagger).$$

	This induces a  norm $\norm{A}=\sqrt{ \inner{A}{A} }$ on the Hilbert space $M_k$ which will be called the Hilbert-Schmidt norm.

\begin{defi} \label{kuIO}
	Given a linear operator $\Phi$ on $M_k$ we denote by $\Phi^*: M_k \to M_k$ the dual linear operator in the sense of Hilbert-Schmidt, that is, if for all $X,Y$ we get
	$$ \inner{ \Phi (X)}{Y}\, = \inner{X}{\Phi^* (Y)}.$$
\end{defi}

\medskip

	Consider a measure $\mu$ on the Borel sigma-algebra over $M_k$.
For an integrable transformation  $F: M_k \to M_k$:
			$$ \int_{M_k} F(v) \, \dm(v) = \left( \int_{M_k} F(v)_{i,j} \, \dm(v) \right)_{i,j}, $$

		in which $F(v)_{i,j}$ is the entry  $(i,j)$ of the matrix $F(v)$.
		
\begin{defi} \label{kui}

		Given a measure $\mu$ on $M_k$ and  a measurable funtion $L:M_k\to M_k$, we say that $\mu$ is  $L$-square-integrable, if

			$$\int_{M_k} \norm{L(v)}^2 \,\dm(v) < \infty.$$

  For a fixed $L$ we denote by $\mathcal{M}(L)$ the set of  $L$-square-integrable measures. We also denote $\mathcal{P}(L)$ the set of  $L$-square-integrable probabilities.
		\end{defi}
			
			\medskip

		 $\phi_L$ is well defined for $L \in \mathcal{M}(L)$.

		\medskip

		\begin{proposition}\label{fun_int1}
			Given a measurable function $L:M_k \to M_k$ and  a square integrable  measure $\mu$, then,  the dual transformation $\phi_L^{*}$  is given by
				$$\phi^*_L(\rho) = \int_{M_k} L(v)^\dagger \rho L(v) \, \dm(v).$$
		\end{proposition}

		\begin{defi} \label{rew}
			Given a measure  $\mu$ over $M_k$ and a square integrable transformation $L: M_k \to M_k$ we say that $L$ is a {\bf stochastic square integrable transformation} if
			$$ \phi_L^*(\text{Id}) = \int_{M_k} L(v)^\dagger L(v) \, \dm(v) = \text{Id}.	$$
		\end{defi}

\begin{defi} \label{posu} A linear map $\phi: M_k \to M_k$ is called {\bf positive} if takes positive matrices to positive matrices.
	\end{defi}

\begin{defi} \label{posu1} A positive linear map $\phi: M_k \to M_k$ is called {\bf completely positive}, if for any $m$, the linear map $\phi_m=\phi \otimes I_m: M_k \otimes M_m \to M_k \otimes M_m$ is positive, where $I_m$ is the identity operator acting on the matrices in $M_m$.
\end{defi}

		\begin{defi} If $\phi:M_k\to M_k$ is square integrable and satisfies
		\begin{enumerate}
				\item $\phi$ is completely positive;
				\item $\phi$ preserves trace,

				then, we say that $\phi$ is a {\bf  quantum channel}.
			\end{enumerate}
		\end{defi}

		\begin{theorem}\label{quantum_channel_phi}

			 Given $\mu$ and $L$ square integrable then the  associated transformation $\phi_L$  is completely positive. Moreover, if  $\phi_L$ is stochastic, then it preserves trace.

		\end{theorem}

		For the proof see \cite{BKL}.
		
		\medskip
		
		\begin{remark}  $\phi_L^*$ is also completely positive.
We say that $\phi_L$ preserves unity if $\phi_L(I) = I$. In this case, $\phi_L^*$ preserves trace. If $\phi_L^*$ preserves the identity then $\phi_L$ preserves trace.
		\end{remark}

		\begin{defi}[Irreducibility]
			We say that $\phi: M_k \to M_k$ is an {\bf irreducible channel} if one of the equivalent properties is true
			\begin{itemize}
				\item Does not exist $\lambda > 0$ and a projection  $p$ in a proper subspace of $\C^k$, such that, $\phi(p) \le \lambda p$;
				\item For all non null $A\ge 0$, $(\textbf{1}+\phi)^{k-1}(A) > 0$;
				\item For all non null $A\ge 0$ there exists $t_A>0$, such that, $(e^{t_A \phi})(A) > 0$;
				\item For all pairs of non null positive matrices  $A,B\in M_k$ there exists a natural number $n\in \{1,...,k-1\}$, such that, $\tr[B\phi^n(A)] > 0$.
			\end{itemize}
		\end{defi}

		For the proof of the equivalences we refer the reader to \cite{evans1978spectral}, \cite{schrader2000perron} and  \cite{wolf2012quantum}.

		\begin{defi}[Irreducibility] \label{gil}  Given $\mu$, we will say (by abuse of language) that $L$ is irreducible if the associated $\phi_L$ is an irreducible channel.
			\end{defi}

		\begin{theorem}[Spectral radius  of $\phi_L$ and $\phi_L^*$]\label{raio_espec}
			Given a square integrable $L:M_k\to M_k$ assume that  the associated $\phi_L$ is irreducible. Then, the spectral radius $\lambda_L> 0$ of $\phi_L$ and $\phi_L^*$ is the same and the eigenvalue is simple. We denote, respectively, by $\rho_L > 0$ and $\sigma_L > 0$, the eigenmatrices, such that,  $\phi_L(\rho_L) = \lambda_L\rho_L$ and $\phi_L^*(\sigma_L) = \lambda_L\sigma_L$, where $\rho_L$ and $\sigma_L$ are the unique non null eigenmatrices  (up to multiplication by scalar).
		\end{theorem}

	The above theorem is the natural  version of the Perron-Frobenius Theorem for the present setting.

	It is natural to think that $\phi_L$ acts on density states and $\phi_L^*$ acts in
	selfadjoint matrices.

\begin{defi}  Given the measure $\mu$  over $M_k$ we denote by $\mathfrak{L}(\mu)$ the set of all integrable $L$ such  that  the associated  $\phi_L$ is irreducible.

\end{defi}

			\begin{defi}
				Suppose   $L$ is in $\mathfrak{L}(\mu)$. We say that $L$ is {\bf normalized} if $\phi_L$ has spectral radius $1$ and preserves trace. We denote by $\mathfrak{N}(\mu)$ the set of all normalized $L$.
			\end{defi}

		If $L \in  \mathfrak{N}(\mu)$, then, we get from  Theorem \ref{raio_espec} and the fact that  $\phi_L^*(I) = I$,  that $\lambda_L = 1$. That is, there exists $\rho_L$ such that $\phi_L(\rho_L) = \rho_L$ and  $\rho_L$ is the only fixed point. Moreover, the spectral radius is equal  to $1$.
		
		\begin{theorem}[Ergodicity and temporal means]
			Suppose $L \in \mathfrak{N}(\mu)$. Then, for all density matrix $\rho \in M_k$ it is true that
			$$ \lim_{N \to \infty} \frac{1}{N}\sum_{n=1}^{N} \phi_L^n(\rho) = \rho_L,$$
			in which $\rho_L$ is the density matrix associated to $L$.
		\end{theorem}
		\textbf{Proof:} The proof follows from  Theorem \ref{raio_espec} and Corollary  6.3 in \cite{wolf2012quantum}.

		\qed\\

		The above result connects irreducibility and ergodicity (the temporal means have a unique limit).

	\section{Entropy} \label{ent}

A measure $\mu$ over $M_k$ (which plays the role of the {\it a priori} probability) is fixed. Thus, given $L\in \mathfrak{L}(\mu)$, we will associate in a natural way the transformation $\phi_L: M_k \to M_k$.

\begin{defi}
We denote by $\Phi=\Phi_\mu$ the set of all $L$ such that the associated
$\phi_L: M_k \to M_k$ is irreducible and stochastic.
\end{defi}

		Suppose  $L$ is irreducible and stochastic.

Given $L$ consider  the density matrix $\rho_L$, which is invariant for $\phi_L$ (see Theorem \ref{raio_espec}).

 \begin{defi}
 We define entropy for $L$ (or, for $\phi_L$) by the expression (when finite) :

		$$h(L)= h_\mu(L) := - \int_{M_k \times M_k} \tr(L(v)\rho_L L(v)^{\dagger}) P_L(v,w)\log P_L(v,w) \,\,\dm(v)\dm(w),$$

		in which
		$$ P_L(v,w) := \frac{\tr(L(w)L(v)\rho_L L(v)^{\dagger}L(w)^{\dagger})}{\tr(L(v)\rho_L L(v)^{\dagger})}.$$
		\end{defi}

		This definition is a generalization of the analogous concept presented on the papers \cite{baraviera2010thermodynamic}, \cite{BLLT2} and \cite{BLLT1}.

The example presented in \cite{BKL} shows that the above definition of entropy is indeed a natural generalization of the classical one in Ergodic Theory. The same example will be considered again when analyzing Lyapunov exponents (see section \ref{exam}).

		\section{ Process $X_n$, $n\in \mathbb{N}$, taking values on $ P(\C^k)$ } \label{pro1}

	Consider a fixed measure $\mu$ on $M_k$ and a fixed $L:M_k \to M_k$, such that, $\int_{M_k} \norm{L(v)}^2\,\dm(v)  < \infty$, and, also assume  that $\phi_L$ is irreducible and stochastic.

Note that if, for example, $\mu$ is a probability and the the function $v \to \norm{L(v)}$ is bounded, we get that
$\int_{M_k} \norm{L(v)}^2\,\dm(v)  < \infty$.

		We follow the notation  of \cite{benoist2017invariant} (and, also \cite{BKL})

		Denote by $P(\C^k)$ the projective space on  $\C^k$ with the metric $d(\hat{x},\hat{y})	= (1-\abs{\inner{x}{y}}^2)^{1/2}$, in which $x,y$   are representatives with norm $1$ and $\inner{\cdot}{\cdot}$ is the canonical inner product.

Take $\hat{x} \in P(\C^k)$ and $S\subset P(\C^k)$. For a stochastic  $\phi_L$ we consider the kernel

	\begin{equation} \label{jajo} \Pi_L(\hat{x}, S) = \int_{M_k} \textbf{1}_{S}(L(v)\cdot\hat{x}) \norm{L(v)x}^2 \,\dm(v),\end{equation}
where the norm above is that of Hilbert-Schmidt.

The discrete-time process (described  by the kernel) taking values on $P(\C^k)$ is determined  by  $\mu$ and $L$. If $\nu$ is a probability on the Borel $\sigma$-algebra $\mathcal{B}$ of $P(\C^k)$ define

		\begin{align*}
			 \nu\Pi_L(S) &= \int_{P(\C^k)} \Pi_L(\hat{x}, S) \,\, \dn(\hat{x}) \\
			 	&= \int_{P(\C^k) \times M_k} \textbf{1}_{S}(L(v)\cdot\hat{x}) \norm{L(v)x}^2 \,\, \dn(\hat{x})\dm(v).
		\end{align*}

		$\nu\Pi_L$ is a new probability on $ P(\C^k)$ and $\Pi_L$ is a Markov operator. The definition above of $\nu \to \nu \Pi_L$ is a simple generalization  of the one in \cite{benoist2017invariant}, in which the authors use the $L$ mentioned here, as the identity transformation.\\

		The map  $\nu \to \nu\,\Pi_L$ (acting on probabilities $\nu$) is called the Markov operator  obtained from  $\phi_L$  in the paper \cite{lozinski2003quantum}, in which the {\it a priori} measure $\mu$ is a sum of Dirac probabilities. Here we consider a more general setting.

		\begin{defi}
			We say that the probability $\nu$ over $P(\C^k)$  is invariant for  $\Pi_L$, if $\nu\Pi_L = \nu$.
		\end{defi}

\medskip
		See \cite{BKL} (and also future Proposition \ref{good}) for the existence of   invariant probabilities for $\Pi_L$

\medskip

		\section{ The process $\rho_n$, $n\in \mathbb{N}$, taking values on $D_k$} \label{pro2}

		For a fixed $\mu$ over $M_k$ and $L$ such $\phi_L$ is irreducible and stochastic, one can naturally define a process $(\rho_n)$ on $D_k = \{\rho\in M_k : \tr \rho = 1 \text{ and } \rho \ge 0\}$, which is called \emph{quantum trajectory}.

		We follow the notation  of \cite{benoist2017invariant} (and, also \cite{BKL})

		Given a $\rho_0$ initial state, we get
$$\rho_n = \frac{L(v)\rho_{n-1}L(v)^*}{\tr(L(v)\rho_{n-1}L(v)^*)}$$ with probability $$\tr(L(v)\rho_{n-1}L(v)^*)\dm(v),\,\,\,\,\,n \in \mathbb{N}.$$

We want to relate the invariant probabilities of the previous section with the fixed point $\rho_{inv} = \rho_{inv}^L$  of $\phi_L$.

\medskip

		First, denote $\Omega := M_k^\N$, and for $\omega=(\omega_i)_{i\in\N}$, take $\phi_n(\omega) = (\omega_1,...,\omega_n)$.

		We denote $\pi_n$ the projection of $\omega$ in its first $n$ coordinates.
		
		We also denote by  $\mathcal{M}$ the Borel sigma algebra  $M_k$. For all, $n\in\N$, consider $\mathcal{O}_n$ the sigma algebra on $\Omega$ generated by the cylinder sets of size $n$, that is, $\mathcal{O}_n := \pi^{-1}_n(\mathcal{M}^{n})$. We equip $\Omega$ with the smaller sigma algebra $\mathcal{O}$ which contains all $\mathcal{O}_n$, $n \in \mathbb{N}$.

		 Denote $\mathcal{J}_n := \mathcal{B}\otimes\mathcal{O}_n$ and $\mathcal{J} := \mathcal{B}\otimes\mathcal{O}$. In this way, $(P(C^k)\times\Omega, \mathcal{J})$ is an integrable space. By abuse of language, we  consider $V_i: \Omega \to M_k$ as a random variable  $V_i(\omega) = \omega_i$. We also introduce another random variable
		
		 $$W_n := L(V_n)...L(V_1), \,\text{\,in which}\,\, W_n(\omega) = L(\omega_n)...L(\omega_1).$$

We point out that the symbol $\otimes$ does not represent the tensor product.

		For a given a probability  $\nu$ on $P(\C^k)$, we define for  $S\in\mathcal{B}$ and $O_n\in\mathcal{O}_n$, another probability

		$$\Pr_\nu(S\times O_n):=\int_{S \times O_n} \norm{W_n(\omega)x}^2\,\,\dn(\hat{x})\dm^{\otimes n}(\omega).$$

	Denote $\E_\nu$ the expected value with respect to $\Pr_\nu$. Now observe that for a $\nu$ probability on $P(\C^k)$, if $\pi_{X_0}: \mathbb{C}^k - \{0\} \to P (\mathbb{C}^k )$ is the orthogonal projection on subspace generated by $X_0$ on $\C^k$, we have

	$$ \rho_{\nu} := \E_\nu(\pi_{X_0}) = \int_{P(\C^k)} \pi_{x_0}\,\, \dn(x_0). $$

	We call $\rho_{\nu}$ barycenter of $\nu$, and it is easy to see	that $\rho_\nu \in D_k$.\\

\medskip

	\begin{proposition}\label{barycenter} If $\nu$ is invariant for $\Pi_L$, then
		$$\rho_\nu = \E_\nu(\pi_{\hat{X}_0}) = \E_\nu(\pi_{\hat{X}_1}) = \phi_L(\rho_\nu).$$

		Therefore,  for an irreducible $L$, every invariant measure $\nu$ for $\Pi_L$ has the same barycenter.
	\end{proposition}

	We emphasize that, in this way, we can recover $\rho_{inv}$, the fixed point of $\phi_L$, by taking the barycenter of any invariant probability (the quantum channel $\phi_L$ admits only one fixed point). That is, for any invariant probability   $\nu$ for $\Pi_L$, we get that $\rho_{\nu}=\rho_{inv}$.

	Note that the previous process can be seen as $\rho_n:\Omega \to D_k$, such that, $\rho_0(\hat{x},\omega) = \rho_\nu$ and

		$$\rho_n(\omega) =\frac{W_n(\omega)\rho_0 W_n(\omega)^*}{\tr(W_n(\omega)\rho_0 W_n(\omega)^*)}.$$

Using an invariant $\rho$ we can define a Stationary Stochastic Process taking values on $M_k$. That is, we will define a probability $\Pr$ over  $\Omega = (M_k)^\mathbb{N}.$

Take $O_n \in \mathcal{O}_n$ and define

	$$\Pr^\rho (O_n) = \int_{O_n} \tr(W_n(\omega)\rho W_n(\omega)^*) \,\,\dm^{\otimes n}(\omega).$$

	The probability $\Pr$ on $\Omega$ defines a Stationary Stochastic Process.

		 \section{Irreducibility, the $\Phi$-Erg property and the purification condition} \label{gen}

In this section,  for the benefit of the reader, we recall the main results from \cite{benoist2017invariant}.
We will use in this section the notation of \cite{benoist2017invariant}.

 \begin{defi}
      Given $L:M_k\to M_k$, $\mu$ on $M_k$ and $E$ subspace of $\C^k$,  we say that $E$ is $(L,\mu)$-invariant, if $L(v)E \subset E$, for all $v\in\supp\mu$.
    \end{defi}

    \begin{defi}
      Given $L: M_k\to M_k$, $\mu$ on $M_k$, we say that $L$ {\bf  is $\Phi$-Erg for $\mu$}, if there exists a unique minimal non-trivial space $E$, such that, $E$ is $(L, \mu)$-invariant.
    \end{defi}

\medskip

In \cite{schrader2000perron}  it is shown that if the above  space $E$ is equal to $\C^k$, then $L$ is {\bf irreducible} for $\mu$ (or, $\mu$-irreducible) in the sense of Definition \ref{gil}.

\medskip

	The relation of $ \Pr^\rho$ and $\Pr_\nu$ (described in the last sections) is described in the next result.
	
	\begin{proposition}  \label{est} The marginal of  $\Pr_\nu$ on $\mathcal{O}$ is $ \Pr^{\rho_\nu}$. If $\Phi$-Erg is true, then for any two $\Pi$-invariant probabilities  $\nu_a$ and $\nu_b$, we get $ \Pr^{\rho_{\nu_a}}= \Pr^{\rho_{\nu_b}}.$
	
	\end{proposition}

The proof of the result above, when $L$ is the identity, was obtained in Proposition 2.1 in \cite{benoist2017invariant}. The proof for the case of a general $L$ is analogous.

\medskip

Given two operators $A$ and $B$ we say that  $A \,\propto \,B$, if there exists $\beta \in \mathbb{C}$, such that, $A = \beta\, B.$

    \begin{defi} \label{kel}
      Given $L: M_k\to M_k$, $\mu$ on $M_k$, we say that the pair $(L,\mu)$ satisfies the {\bf purification condition}, if an orthogonal projector $\pi$, such that, for any $n \in \mathbb{N}$
      $$ \pi L(V_1)^*...L(V_n)^*\,  L(V_n)...L(V_1)\,\pi\,\, \propto\,\, \pi,$$
      for $\mu^{\otimes n}$-almost   all $ (v_1,v_2,..,v_n)$, is necessarily of rank one.
\end{defi}

\medskip

Following \cite{benoist2017invariant}, we denote $\mathbb{P}^{ch}= \mathbb{P}^{ \frac{1}{k} \, Id}.$

We denote by $Y_n$, $n \in \mathbb{N}$, the matrix-valued random variable
$$ Y_n \, = \frac{W^*_n \,W_n  }{ \text{tr}\,\,(W^*_n \,W_n)}, \,\,\text{if\,\,tr } (W^*_n \,W_n) \neq 0,$$
that we extend the definition arbitrarily when tr  $(W^*_n \,W_n)$=0 .

\medskip

The next two propositions are of fundamental importance for the theory and they were  proved in Proposition 2.2 in \cite{benoist2017invariant} (the same proof applies here).

\begin{proposition}
For any probability $\nu$ over $P(\mathbb{C}^k)$  the stochastic process $Y_n$, $n \in \mathbb{N}$, is a martingale with respect to the sequence of sigma-algebras $\mathcal{O}_n$, $n \in \mathbb{N}$. Therefore, there exists a random variable $Y_\infty$ which is the  almost sure limit of $Y_n$ for the probability $\mathcal{P}_\nu$ and also in the $L^1$ norm.

\end{proposition}

\begin{proposition}

For any probability $\nu$ over $P(\mathbb{C}^k)$ and $\rho\in \mathcal{D}_k$
$$ \frac{d\, \mathbb{P}^\rho  }{d\, \mathbb{P}^{ch}} = k\, \,\text{tr}\,\, (\rho\, Y_\infty).$$

Moreover, $\mu$ and $L$ satisfy the purification condition, if and only if, $Y_\infty$ is $\mathbb{P}_\nu$-a.s a rank one projection for any probability $\nu$ over $P(\mathbb{C}^k)$.

\end{proposition}

\begin{proposition} \label{good} If the pair $(L,\mu)$ satisfies the $\phi$-Erg and the purification condition, then, the Markov kernel $\Pi$ admits a unique invariant probability.

\end{proposition}

$x_1\wedge x_2 \wedge ... \wedge x_n$, with $x_j\in \mathbb{C}^k$, denotes the classical
wedge product (an alternate form on $\mathbb{C}^k$).

One can consider an inner product
$$    \langle\, r_1\wedge r_2 \wedge ... \wedge r_n\,\,, \, \,s_1\wedge s_2 \wedge ... \wedge s_n \,\ \rangle   = \,\,\text{det}\,\, ( r_i s_j ) )_{i,j=1,2,..., n} ,
$$
and, the associated norm $ | x_1\wedge x_2 \wedge ... \wedge x_n\,|.$

Given an operator $X: \mathbb{C}^k \to \mathbb{C}^k $ we define
$ \bigwedge^n X \, : \wedge^n \mathbb{C}^k \to  \wedge^n \mathbb{C}^k $ by
$\bigwedge^n X \,( x_1\wedge x_2 \wedge ... \wedge x_n)\,= X( x_1)\wedge X(x_2) \wedge ... \wedge X(x_n).$

\begin{proposition} \label{goo}
Assume the pair $(L,\mu)$ satisfies the purification condition, then, there are two constants $C>0$ and $\beta<1$, such that, for each $n$
$$\int_{M_k^n}\, |\,\wedge^2 \,(\,L(v_n) ... L(v_1) \,)\, | \,d\, \mu^{\otimes \,n} (v_1,v_2,...,v_n)\, = \mathbb{E}^{ch} \,(k\, \frac{| \, \bigwedge^2 \,W_n\,| }{\text{tr}\,\, (\,W_n^*\, W_n \,)} \,)\leq C\, \beta^n.$$

\end{proposition}

The proofs of the two propositions above are similar to the corresponding ones in \cite{benoist2017invariant}.

\medskip

    Consider $\Bcal(M_k) = \{ L: M_k \to M_k \,\vert\, \text{L is continuous and bounded} \}$ where $\norm{L} = \sup_{v\in M_k} \norm{L(v)}$.
    \begin{defi}
      For a fixed measure $\mu$ over $M_k$, define

      $$\Bcal_\mu(M_k) = \{ L\in \Bcal \,\vert\, \text{L is $\mu$-irreducible}\},$$

      and

      $$\Bcal^\Phi_\mu(M_k) = \{ L\in \Bcal \,\vert\, \text{L is $\Phi$-Erg for $\mu$} \}.$$
    \end{defi}

    \begin{proposition}\label{phi_erg_dense_set}
      Given $\mu$ over $M_k$ with $\#\supp\mu > 1$, $\Bcal^\Phi_\mu(M_k)$ is open and dense on $\Bcal(M_k)$.
    \end{proposition}
    \textbf{Proof:} See  \cite{BKL}.

    \qed

\medskip

In section \ref{purgen} we will prove:
\begin{proposition}\label{phi_erg_dense_set}
      Given $\mu$ over $M_k$ with $\#\supp\mu > 1$, the set of $L$ satisfying the purification condition is generic in $\Bcal(M_k)$.
    \end{proposition}

	\section{Lyapunov exponents for quantum channels} \label{Ly}

	In this section, we will consider a discrete-time process taking values on $M_k$.
	
	Take $\mu$ over $M_k$ and  $L: M_k \to M_k$ in such way that the associated channel $\Phi$ defines a $\Phi$-Erg stochastic map. We assume in this section that $\rho \in D_k$ is  such that $\Phi(\rho) = \rho$.
	Such $\rho$ plays the role of the initial vector of probability (in the analogy with the theory of Markov Chains).
	
	We follow the notation of \cite{benoist2017invariant}.

	Take $\Omega = M_k^{\N}$, and for $n\in \N$ let $\mathcal{O}_n$ be the $\sigma$-algebra on $\Omega$ generated by the $n$-cylinder sets (as in Section \ref{pro2}).

	An element on $\Omega$ is denoted by $(\omega_1,\omega_2,...,\omega_n,...).$
Following section \ref{pro2} we denote
	$W_n(\omega) = L(\omega_n)...L(\omega_1).$

	Taking $O_n \in \mathcal{O}_n$ we define

	$$\Pr(O_n) = \int_{O_n} \tr(W_n(\omega)\rho W_n(\omega)^*) \,\,\dm^{\otimes n}(\omega).$$

	If $\mathcal{O}$ is the smallest $\sigma$-algebra of $\Omega$ that contains all $\mathcal{O}_n$ we can extend the action of $\Pr$  to this $\sigma$-algebra.
	
	The probability $\Pr$ on $\Omega$ defines a Stationary Stochastic Process.

	\begin{theorem}
		$(\Omega, \Pr, \theta)$ is ergodic where $\theta$ is the shift map.
	\end{theorem}

	The above theorem has been proved in Lemma 4.2 in \cite{benoist2017invariant}.
	
	\medskip

	\begin{theorem} \label{pit}
		Suppose the pair $(L,\mu)$ satisfies irreducibility, the $\phi$-Erg and the purification condition. Assume also that $\int |L(v)|^2 \log |L(v)|\, \ d\mu(v)<\infty$, then, there exists numbers
		$$\infty > \gamma_1\geq \gamma_2\geq ...\geq \gamma_k\geq -\infty,$$
		such that, for any probability $\nu$ over $P(\mathbb{C}^k$) and any $p \in \{1,2,..,k\}$
		$$ \lim_{n \to \infty} \frac{1}{n} \log | \bigwedge^p \, W_n\,| =\sum_{j=1}^p \gamma_j,$$
		$\mathbb{P}_{\nu}$-a.s.
		\end{theorem}

	The above theorem was proved in \cite{benoist2017invariant} and the same proof works in our current setting. We point out that a key element in this proof (see (35) in \cite{benoist2017invariant}) is the fact that if $(L,\mu)$ is $\phi$-Erg and
	irreducible, then, $\rho_{inv}>0$, and for any $\rho\in \mathcal{D}_k$
	$$ \mathbb{P}^\rho << \mathbb{P}^{\rho_{inv}}.$$

	Proposition \ref{est} is also used in the proof (corresponds to proposition 2.1 in \cite{benoist2017invariant}).
	
	The numbers
	$$\gamma_1\geq \gamma_2\geq ...\geq \gamma_k,$$
	are called the {\bf  Lyapunov exponents}.
	\begin{theorem}
	Suppose that $L$ is generic for $\mu$ and $\int |L(v)|^2 \log |L(v)|\, \ d\mu(v)<\infty$, then, the Lyapunov exponents $$\infty > \gamma_1\geq \gamma_2\geq ...\geq \gamma_k\geq -\infty$$ are well defined.
	
	(a) $\gamma_2 -\gamma_1<0$, where $\gamma_2-\gamma_1$ is the limit
	$$  \lim_{n \to \infty } \frac{1}{n}  \log \frac{ |\,\bigwedge^2 W_n\, |}{|\,W_n\,|^2 },$$
	whenever $\gamma_1 = -\infty$

	\medskip

	(b) For  $\mathbb{P}_\nu$-almost sure $x$ we have that
	$$  \lim_{n \to \infty } \frac{1}{n} ( \log \,|\,W_n\, (x)\,| -  \log \,|\,W_n\,| )=0.$$
	
\end{theorem}

\textbf{proof:}

	The proof of (a) is similar to the one in \cite{benoist2017invariant}.
	
	\medskip

Proof of  (b): We start with

%Since $M_n \in \mathcal{D}_k$, we have $\tr(M_n)=1$ and

$$\frac{\norm{W_n x}}{\norm{W_n}}=\frac{[ \tr (W_n \pi_{x} W_n^{\dagger})\ ]^{1/2}}{[ \tr  (W_n W_n^{\dagger})\ ]^{1/2}} = \left[ \frac{ \tr (W_n^{\dagger}W_n \pi_{x})}{ \tr  (W_n W_n^{\dagger})}  \right]^{1/2}$$

$$ = \left[  \tr \left( \dfrac{W_n^{\dagger} W_n}{\tr  (W_n W_n^{\dagger})}\ \pi_{x}\right)  \  \right]^{1/2} = \left[  \tr \left( M_n  \pi_{x}\right)   \right]^{1/2}  . $$

Those calculations are valid for all $\omega$ such that $W_nW_n^{\dagger}(\omega) \neq 0$. Since $\Pr_{\nu}(W_n = 0) = 0$, no extra work is required. Now we use proposition 2.2 in \cite{benoist2017invariant} which says that $M_n$ converges $\Pr_{\nu}-a.s.$ and in $L^1$ norm to a $\mathcal{O}$-measurable random variable $M_{\infty}$. By continuity of the trace and square root, we have
$$ \lim_{n \to \infty } \left[  \tr \left( M_n  \pi_{x}\right)   \right]^{1/2} = \left[  \tr \left( M_{\infty}  \pi_{x}\right)   \right]^{1/2}, \text{ for } \Pr_{\nu}-\text{a.s. } x \in P(\C ^k). $$

The proof is similar to the one in \cite{benoist2017invariant}.

\qed

\medskip

\section{The main example - Lyapunov exponents and entropy} \label{exam}

Now we will present an example  (see example 8.5 in \cite{BKL})  in which we can  estimate the Lyapunov exponents and show a relation with entropy.
We will call it the {\it  Markov model in quantum information}.
\smallskip

A concept of entropy for Channels was introduced in \cite{baraviera2010thermodynamic} and this concept is quite natural because extends the concept of entropy for Markov Chains. Given a stochastic matrix $P = (p_{ij})$ we can associate to it an invariant Markov probability measure. 

We will consider a particular case (see example 8.5 in \cite{BKL})  where $P$ is a $2\times 2$ matrix.
Our main result in this section is for the first  Lyapunov exponent $\gamma_1$:

\begin{equation} \label{saca1} \gamma_1 = \frac{1}{2}\sum_{i,j \in \{0,1\}}\pi_j p_{ij}\log(p_{ij}) = -\frac{1}{2} h,
		\end{equation}
in which $h$ is the entropy of the Markov invariant measure associated to the $2\times 2$ matrix$P$.

		Let $V_{ij} = \sqrt{p_{ij}}\ket{i}\bra{j}$ where $P = (p_{ij})$ is a irreducible (in the classical sense for a Markov chain) $k$ by $k$ column stochastic matrix, $\mu = \sum_{ij} \delta_{V_{ij}}$ and $L=I$. In this case, we get that
				$$V_{ij}^* V_{ij} = p_{ij}\ket{j}\bra{j},$$
		is a diagonal matrix and, if $A=(a_{ij})$,

				$$V_{ij}^* A V_{ij} = p_{ij}a_{ii} \ket{j}\bra{j}.$$

		Therefore, when $\omega=(V_{i_n j_n})$, we have

				$$W_2(\omega)^*W_2(\omega) = p_{i_1j_1}p_{i_2j_2} \ket{j_2}\bra{j_2}\delta_{i_2 j_1},$$
		in which $\delta_{ij} = 0$ if $i\ne j$ and $1$ if $i=j$. By induction,
				$$W_n(\omega)^*W_n(\omega) = \left(\prod_{k=1}^{n} p_{i_kj_k}\right) \left(\prod_{k=1}^{n-1} \delta_{i_{k+1}j_k}\right) \ket{j_n}\bra{j_n}.$$

		Thus, $W_n(\omega)^*W_n(\omega)$ is $0$ or a diagonal matrix with a unique entry different from $0$. This entry is exactly $\left(\prod_{k=1}^{n} p_{i_kj_k}\right)$ which implies that
			$$\norm{W_n(\omega)^*W_n(\omega)} = \left(\prod_{k=1}^{n} p_{i_kj_k}\right) = p_{11}^{X_{11,n}(\omega)}p_{12}^{X_{12,n}(\omega)}p_{21}^{X_{21,n}(\omega)}p_{22}^{X_{22,n}(\omega)},$$
		in which

		$$X_{ij,n}(\omega) = \sum_{k=0}^{n-1} \1_{[V_{ij}]} \circ \theta^k(\omega), $$
with $\1_{[V_{ij}]}$ being the characteristic function of cylinder $[V_{ij}]$.

		Note that under the  ergodicity hypothesis we get the property: for any $i,j$ and $\Pr$-almost sure $\omega$, we have that
	$$	\lim_{n \to \infty} \frac{1}{n} \, \sum_{k=0}^{n-1} \1_{[V_{ij}]} \circ \theta^k(\omega)= \Pr([V_{ij}])  .$$
		
		It follows from the arguments of example 8.5 in \cite{BKL} that the pair $(L,\mu)$ satisfies the purification condition, the $\phi$-Erg conditions, and also irreducibility.

		 Remember that for given matrix $A \in M_k(\C)$, $a_1(A) \ge a_2(A) \ge ... \ge a_k(A)$ are the singular values of $A$, i.e., the square roots of eigenvalues of $A^*A$, labeled in decreasing order. From Lemma III.5.3 in \cite{lacroix} we have $$\norm{\bigwedge^p W_n(\omega)} = a_1(W_n)...a_p(W_n).$$
		
		Therefore,  $\norm{\bigwedge^1 W_n(\omega)} = a_1(W_n) = \norm{W_n^*W_n}^{\frac{1}{2}}$.
		
		Following Proposition \ref{pit}  (which corresponds to Proposition 4.3 in \cite{benoist2017invariant}) we can obtain the greater Lyapunov exponent $\gamma_1$ taking the limit

		\begin{align*}
		\gamma_1 :&= \lim_{n \to \infty} \frac{1}{n} \norm{\bigwedge^1 W_n(\omega)}\\
			&= \lim_{n \to \infty} \frac{1}{n} \log(\norm{W_n(\omega)^*W_n(\omega)}^{\frac{1}{2}})\\
			&= \frac{1}{2}\sum_{ij} \lim_{n \to \infty} \frac{X_{ij,n}(\omega)}{n} \log(p_{ij})\\
			&= \sum_{ij} \frac{\Pr([V_{ij}])}{2}\log(p_{ij}).
		\end{align*}

We know from  \cite{BKL} that if $\Phi$ is the channel defined for such $\mu$ and $L$, then, $\Phi$ is $\Phi$-Erg. Assume, $\pi = (\pi_1, \pi_2,...,\pi_k)$ is the unique invariant probability vector for the stochastic matrix $P$.
Then, the unique $\rho$, such that, $\Phi(\rho) = \rho$, is exactly
	the diagonal matrix 	$\rho$ with entries $\pi_1, \pi_2,...,\pi_k$. We also know, in this case, that the entropy (see example 8.5 in \cite{BKL}) for a channel as defined in \cite{BKL} is equal to the classical Shannon-Kolmogorov entropy of the stationary Markov Process associated to the  column stochastic matrix $P = (p_{ij})$ (see formula in \cite{PY}).

Now, we can estimate
\begin{align*}
			\Pr([V_{ij}]) &= \int_{[V_{ij}]} \tr(v\rho v^*) \,\dm(v)
			= \tr(V_{ij}^*V_{ij}\rho)
			= \tr(p_{ij}\ket{j}\bra{j}\rho)
			= p_{ij}\bra{j}\rho\ket{j}\\
			&= p_{ij}\pi_j.
		\end{align*}

		Therefore,

\begin{equation} \label{saca} \gamma_1 = \frac{1}{2}\sum_{i,j \in \{0,1\}}\pi_j p_{ij}\log(p_{ij}) = -\frac{1}{2} h,
		\end{equation}
in which $h$ is the entropy of the Markov invariant measure associated to the matrix $P$.

		The value $\frac{1}{2}$ which multiplies the entropy on the above expression is due to the fact that we considered the norm $\norm{A}=\inner{A}{A}^{1/2}$.

Now we estimate the second Lyapunov exponent $\gamma_2$.

		We showed that
		$W_n(\omega)^*W_n(\omega) = \left(\prod_{k=1}^{n} p_{i_kj_k}\right) \left(\prod_{k=1}^{n-1} \delta_{i_{k+1}j_k}\right) \ket{j_n}\bra{j_n}$, which implies that the second eigenvalue is $0$ and therefore $a_2(W_n(\omega))=0$.
		
		Now, we can get $\gamma_2$, in fact,

		$$\gamma_1 + \gamma_2 = \lim_n\frac{1}{n}\log(a_1(W_n(\omega))a_2(W_n(\omega)))=\lim_n\frac{1}{n}\log(0),$$

		which implies that $\gamma_2 = -\infty$.

\medskip

\section{The purification condition is generic} \label{purgen}

\medskip
The measure $\mu$ is fixed from now on.

Our main goal in this section is to show:

\begin{proposition}\label{phi_erg_dense_set1}
      Given $\mu$ over $M_k$ with $\#\supp\mu > 1$, the set of $L$ satisfying the purification condition is generic in $\Bcal(M_k)$.
    \end{proposition}

    This will follow from Lemma \ref{mainle}.

  \begin{defi}
    We say that the projection $\pi$ $n$-purifies $L:M_k\to M_k$, in which $\rank\,\pi\ge 2$, if there exists $E \in \mathcal{O}_n$, with $\mu^{\otimes n}(E) > 0$, such that,

    $$ \pi W_n(\omega)W_n(\omega)\pi \not\propto \pi ,$$

    for all $\omega \in E$.
  \end{defi}

In order to prove that a certain $L$ satisfies the purification condition, we have to consider all possible projections $\pi$ (see definition \ref{kel}).

\medskip

  Observe that if $Q$ is a unitary matrix, $\pi$ has rank great or equal to $ 2$ and $n$-purifies $L$ for $E \in \mathcal{O}_n$, then

  $$ Q\pi Q^*Q W_n(\omega)^*W_n(\omega)Q^* Q\pi Q^* \not\propto Q \pi Q^*. $$

  Besides that, $W_n(\omega) = L(\omega_n)...L(\omega_1)$, so, as $Q^*Q = I$, we have

  $$QW_n(\omega)^*W_n(\omega) Q^* =$$
  $$QL(\omega_1)^*Q^*Q...Q^*QL(\omega_n)^*Q^*QL(\omega_n)Q^*Q...Q^*QL(\omega_1)Q^*Q.$$

  From this follows:

  \begin{proposition}
    If $L_Q(v) := QL(v)Q^*$, then for a projection $\pi$, such that, $\rank\,\pi \ge 2$ and an unitary matrix $Q$, it's true that

    $$\pi \npurifies L \iff Q\pi Q^* \npurifies L_Q.$$
  \end{proposition}

  \begin{defi}
    For an orthogonal projection $\pi$ and $n\in \N$ we define

    $$ Pur_\pi^n = \{L \in \Bcal(M_k) \, \vert \, \pi \npurifies L\}. $$
  \end{defi}

  Note that if
  $$Pur = \{L \in \Bcal(M_k) \, \vert \, \Phi_L \text{ satisfies (Pur) condition } \},$$
  and if we denote
  $$P_2 = \{\pi \text{ orthogonal projection }\,\vert\, \rank\pi \ge 2\},$$
it follows that
    $$Pur = \bigcap_{\pi \in P_2}\,\bigcup_{n\in\N} Pur_\pi^n.$$

  \begin{proposition}\label{purn_is_open}
    For any $\pi\in P_2$ and $n\in\N$, $Pur_\pi^n$ is open.
  \end{proposition}
  \textbf{proof:}
    Take $\pi \in P_2$ with $\rank \pi = l$, $Q$ an unitary matrix that diagonalizes $\pi$. Suppose that

    \[
    \tilde{\pi} := Q\pi Q^* = \begin{bmatrix}
        1 & 0 & \dots & 0 & 0 & \dots & 0\\
        0 & 1 & \dots & 0 & 0 & \dots & 0\\
        \vdots & \vdots & \ddots & \vdots & \vdots & \ddots & \vdots\\
        0 & 0 & \dots & 1 & 0 & \dots & 0\\
        0 & 0 & \dots & 0 & 0 & \dots & 0\\
        \vdots & \vdots & \ddots & \vdots & \vdots & \ddots & \vdots\\
        0 & 0 & \dots & \dots & 0 & \dots & 0
        \end{bmatrix}.
    \]

  So, if $L \in Pur_\pi^n$, we know that $\tilde{\pi} \npurifies L_Q$.

  If $$ (s^L_{ij}(\omega)):=QW_n(\omega)^*W_n(\omega)Q^* $$
we know that $s^L_{ij}$ are continuous functions. Moreover, there exists $\omega_0\in \supp\mu^{\otimes n}$ such that at least one of following conditions occurs:

  1. There exists $i\ne j$, such that, $s^L_{ij}(\omega_0) \ne 0$, or there exists $j=i>$ rank $\tilde{\pi}$, such that, $s^L_{i i}(\omega_0) \ne 0$.

In this case, we define the matrix $s^L(\omega) := s^L_{ij}(\omega)$; \\

  2. There exists $i \ne j$, such that, $s^L_{ii}(\omega_0) - s^L_{jj}(\omega_0) \ne 0$.

  In this case, we define $s^L(\omega) := s^L_{ii}(\omega) - s^L_{jj}(\omega)$.\\

  It is clear that for $F\in \Bcal(M_k)$ with $\norm{L - F} \le \varepsilon$, we have, for $\varepsilon$ small enough, that $s^F(\omega_0) \ne 0$, because $s^F(\omega_0)$ has a continuous dependence on F. Furthermore, $s^F$ is continuous, so there is a open set $B$ with $\omega_0\in B$, such that, $s^F(\omega)\ne 0$, for $\omega \in B$. Moreover, $\omega_0 \in \supp\mu^{\otimes n}$, which implies that $\mu^{\otimes n}(B) > 0$, and therefore $F \in Pur_\pi^n$.

 \qed

  \medskip
 We point out that $s^L(w_0) \neq 0$ is the good condition for purification.
  \medskip

  \begin{proposition}
    For any $\pi \in P_2$, $Pur_\pi^1$ is dense.
  \end{proposition}
  \textbf{proof:} Take $L \notin Pur_\pi^1$ and $Q$ unitary matrix that diagonalizes $\pi$ as above. If $L \notin Pur_\pi^1$, defining $(s^L_{ij})$ as in the previous proposition, we know that $s_{11}(v)=s_{22}(v)$, for almost every $v$. If $D = \ket{1}\bra{1}$, and $\varepsilon > 0$, we consider $L_Q^\varepsilon(v) = L_Q(v) + \varepsilon D$. So, we have

  \begin{align*}
    L_Q^\varepsilon(v)^*L_Q^\varepsilon(v) &= (L_Q(v) + \varepsilon D)^*(L_Q(v) + \varepsilon D)\\
    &= L_Q(v)^*L_Q(v) + \varepsilon L_Q^*(v)D + \varepsilon DL_Q(v) + \varepsilon^2 D.
  \end{align*}

  This perturbation changes the values $s^L_{11}(v)$ but not $s^L_{22}(v)$, which implies that $Q\pi Q^*$ $\npurifies$  $L_Q^\varepsilon$. Now, we only need to take $L^\varepsilon(v) = L(v) + \varepsilon Q^*DQ$ and, as $\norm{L - L^\varepsilon}$ is small, it will follow the density property at once.

  \qed

  \begin{defi}
    For any $\pi \in P_2$, we define $Pur_\pi = \bigcup_{n\in \N}Pur_\pi^n$.
  \end{defi}

Note that  by the two propositions above given a fixed $\pi$ the set  $Pur_\pi$ is open and dense on $\Bcal(M_k)$.  Note also that for the purification condition it is necessary  to consider all possible projections $\pi$ (see definition \ref{kel}).

  \begin{lemma}
     If $\pi_1, \pi_2 \in P_2$ has the same rank, $E_i := \im \pi_i$, and $\{ x_i \}$ is an orthonormal basis for $E_i$, then, $\{ \pi_2 x_i \}$
      is a basis for $E_2$ if $\pi_1$ and $\pi_2$ are close enough.
  \end{lemma}
  \textbf{proof:}  The proof will be by contradiction. Suppose $\rank \pi_1 = l$ and $y_i := \pi_2 x_i$ are linearly dependent. So, the dimension generated by $\{y_i\}$ is at most $l-1$. Then, there  exists a vector  $\hat{y}\in E_2$, which has norm $1$  and it is orthogonal to the subspace generated by $\{y_i\}$. Therefore,  $\inner{\hat{y}}{y_i} = 0$, for all $i$. This implies that $\inner{\hat{y}}{\pi_2x_i} = \inner{\pi_2\hat{y}}{x_i} = \inner{\hat{y}}{x_i} = 0$ and, moreover, $\pi_1\hat{y} = 0$. Finally, we get  $\norm{\pi_1 - \pi_2} \ge \norm{\pi_1\hat{y} - \pi_2\hat{y}} = \norm{\hat{y}} = 1$.

 If we assume that $\norm{\pi_1 - \pi_2}<1$, we are done.\\

  \qed

  Observe that, for $i\ne j$ and $\norm{\pi_2 - \pi_1} < \varepsilon$,
  $$\abs{\inner{y_i}{y_j}} = \abs{\inner{y_i}{x_j}} = \abs{\inner{y_i-x_i}{x_j}}
      \le \norm{y_i - x_i}\norm{x_j} = \norm{y_i - x_i}$$
 \begin{equation} \label{we} \le \norm{\pi_2 - \pi_1} < \varepsilon
 \end{equation}

  The set of $y_i$ is not an orthonormal basis.

  We would like to get an orthonormal basis close to the orthogonal basis $x_1,..,x_n$.
  Our aim is to prove corollary \ref{rer} which claims  that, given $\varepsilon$ there exists an orthonormal basis  $(u_i)$ for $E_2$ with $\norm{u_i - x_i} < C\varepsilon$, for some constant $C>0$. In this direction, we will perform a Gram-Schmidt normalization procedure.

  Denote $u_1 := \dfrac{y_1}{\norm{y_1}}$, $N_i := y_i - \sum_{j=1}^{i-1} \inner{y_i}{u_j}u_j$  and $u_i :=\dfrac{N_i}{\norm{N_i}}$, for $i > 1$. Then, we have

  \begin{align*}
    \norm{u_1 - x_1} &= \biggnorm{
      \norm{\pi_2x_1}^{-1}\pi_2x_1 -
      \norm{\pi_2x_1}^{-1}x_1 +
      \norm{\pi_2x_1}^{-1}x_1
      - x_1
    } \\
    &\le \norm{\pi_2x_1}^{-1}\norm{\pi_2x_1 - \pi_1x_1} +
         \norm{x_1}\abs{\norm{\pi_2x_1}^{-1} - 1}\\
    &\le \norm{\pi_2x_1}^{-1}\varepsilon +
         \abs{\norm{\pi_2x_1}^{-1} - 1}\\
    &<   \varepsilon(1-\varepsilon)^{-1} + (1-\varepsilon)^{-1}\varepsilon\\
    &<   4\varepsilon,
  \end{align*}

  and
    $$\abs{\norm{\pi_2x_i} - \norm{x_i}} \le \norm{\pi_2x_i - \pi_1x_i} < \varepsilon \implies 1 - \varepsilon < \norm{\pi_2x_i} < 1 + \varepsilon.$$

  Furthermore,
  \begin{align*}
    \norm{y_2 - \inner{y_2}{u_1}u_1} &= \norm{
        y_2 -
        \norm{\pi_2x_i}^{-2}\inner{y_2}{y_1}y_1
    }\\
    &\le \norm{x_2} + (1 - \varepsilon)^{-2}\varepsilon \norm{x_1}\\
    &< 1 + 4\varepsilon.
  \end{align*}

  \begin{proposition}
    For any $i \in \{1,...,n\}$, $j<i$, there is $C_{ij} > 0$, such that, $\abs{\inner{y_i}{u_j}} < C_{ij}\varepsilon$.
  \end{proposition}
 \textbf{proof:}

  Take $N := \min{\norm{N_i}} > 0$.

   Observe that $\abs{\inner{y_2}{u_1}} \le N^{-1}\varepsilon$ (this follows from a similar procedure described in (\ref{we}) and the  Cauchy-Schwartz inequality).

  Then suppose, for all $l < i$, $\abs{\inner{y_l}{u_j}} \le C_{lj}\varepsilon$, and for all $j < l$, with $C_{lj} > 0$. If $j < i$, we have
  \begin{align*}
    \abs{\inner{y_i}{u_j}} &\le N^{-1}\abs{\inner{y_i}{y_j}} + N^{-1}\sum_{k = 1}^{j-1} \abs{\inner{y_j}{u_k}}\\
    &\le N^{-1}\varepsilon + N^{-1}\sum_{k=1}^{j-1} C_{jk}\varepsilon\\
    &=\bigg(1 + \sum_{k=1}^{j-1} C_{jk}\bigg)N^{-1}\varepsilon.
  \end{align*}

 Taking $C_{ij} =\bigg(1 + \sum_{k=1}^{j-1} C_{jk}\bigg)N^{-1}$ the claim follows by induction.

  \qed

  \begin{proposition}
    For all $i$, $\abs{\norm{N_i} - 1} < K\varepsilon$, for some $K>0$.
  \end{proposition}
  \textbf{proof:}
  \begin{align*}
    \norm{N_i} &\le 1 + \norm{y_i - x_i} + \sum_{j=1}^{i-1}\abs{\inner{y_i}{u_j}}\\
    &\le 1 + \varepsilon + \bigg(\sum_{j=1}^{i-1}C_{ij}\sum\bigg)\varepsilon\\
    &= 1 + \bigg(1 + \bigg(\sum_{j=1}^{i-1}C_{ij}\bigg)\bigg) \varepsilon.
  \end{align*}

  Taking $K = \abs{1 + \sum_{j=1}^{i-1}C_{ij}}$  the proof is completed.

  \qed

  \begin{proposition}
    For all $i$, we have $\norm{u_i - x_i} < C_i\varepsilon$.
  \end{proposition}
  \textbf{proof:}
  \begin{align*}
    \norm{u_i - x_i} &\le N^{-1}\norm{N_i - x_i}  + \norm{\norm{N_i}x_i - x_i}\\
    &\le N^{-1}\norm{y_i - x_i} + N^{-1}\sum_{j=0}^{i-1}\abs{\inner{y_i}{u_j}} + \abs{\norm{N_i} - 1}\\
    &\le N^{-1}\varepsilon + N^{-1}\sum_{j=1}^{i-1}C_{ij} + K\varepsilon\\
    &= \bigg(N^{-1} + N^{-1}\sum_{j=1}^{i-1}C_{ij} + K \bigg)\varepsilon.
  \end{align*}

  Define $C_i := N^{-1} + N^{-1}\sum_{j=1}^{i-1}C_{ij} + K$ and the statement has been proved.

  \qed

  \begin{corollary} \label{rer} There exists $C>0$, such that, for all $\varepsilon>0,$
 there exists an orthonormal basis $(u_i)$  for $E_2$ with $\norm{u_i - x_i} < C\varepsilon$.
  \end{corollary}

  If we repeat the process, but now for $E_1^{\perp}$ and $E_2^{\perp}$, we get another constant $C_2$,  and in a similar way, we obtain new vectors $(u_i)$ from the $(x_i)$. These $u_i$ define an orthonormal basis for $\C^k$ with $\norm{u_i - x_i} < C_2\varepsilon$.

Now we define $Q_1, Q_2$ , such that $Q_1x_i = e_i$ and $Q_2u_i = e_i$, in which $(e_i)$ is a canonical basis for $\C^k$. Then, $Q_i\pi_iQ_i^*$ is a diagonalization for $\pi_i$. Observe that
$$\norm{Q_2Q_1^* - I}=\norm{Q_2Q_1^* - Q_2 \, Q_2^*}  < C_3\varepsilon,$$ the map $A \to Q_2Q_1^* A (Q_2Q_1^*)^*$ is continuous and, in addition, this map is close to the identity map.\\

  We know that if $L \in Pur_{\pi_1}^n$, then,  we take $s_1^L$ from the $\pi_1, Q_1$ and $L$, as in  proposition \ref{purn_is_open}. In the same way, there exists $\omega_0 \in \supp\mu$ with $s_1^L(\omega_0) \ne 0$. Observe that $s_1^L(\omega_0)$ depends only on the cordinates of $Q_1W^L_n(\omega_0)^*W^L_n(\omega_0)Q_1^*$. Now, applying $Q_2Q_1^*\,(\cdot) \,(Q_2Q_1^*)^*$ we get $Q_2W^L_n(\omega_0)^*W^L_n(\omega_0)Q_2^*$. Note that this is the same as considering $\pi_2, Q_2, L$ and the associated $s_2^L(\omega_0)$. If $\varepsilon$ is small enough, $s_2^L(\omega_0) \ne 0$ and we can repeat the argument used in proposition \ref{purn_is_open} in order to obtain an open set $B$, such that, $\omega_0 \in B$ and, moreover,  if $\omega\in B$ then $s_2^L(\omega) \ne 0$. Therefore, $L \in Pur_{\pi_2}^n$.\\

  The previous  arguments prove the following lemma.

  \begin{lemma}\label{close_pi_pur}
    If $L \in Pur_{\pi_1}$ and $\pi_1, \pi_2$ are close enough, then $L \in  Pur_{\pi_2}$.
  \end{lemma}

  \begin{lemma} \label{mainle}
    Take {\bf $K_2$ as a countable dense subset of $P_2$}. Then,

    $$Pur = \bigcap_{\pi \in K_2}\,Pur_\pi.$$
  \end{lemma}
  \textbf{proof:} We will use the classical Baire Theorem.

  Suppose that $L \in  \bigcap_{\pi \in K_2}\,Pur_{\pi}$,
  then for lemma \ref{close_pi_pur}, for every $\pi \in K_2$,
  as $L\in Pur_\pi$, there exists $\varepsilon(\pi)$, such that, if
  $\hat{\pi}\in P_2$ and $\norm{\pi - \hat{\pi}} < \varepsilon(\pi)$,
  then $L \in Pur_{\hat{\pi}}$. Furthermore, if we define
  $B(\pi) = \{\hat{\pi}\in P_2 | \norm{\hat{\pi} - \pi} < \varepsilon(\pi)\}$,
  then $\bigcup_{\pi\in K_2} B(\pi)$ covers $P_2$. Therefore,
  for any $\hat{\pi}\in P_2$, there is $\pi \in K_2$, such that, $\hat{\pi}\in B(\pi)$
  and thus $L \in Pur_{\hat{\pi}}$. This implies that $L \in \bigcap_{\pi \in P_2} Pur_{\pi} = Pur$.

  \qed%
\section{Appendix - Relations with previous works and with classical dynamical systems} \label{appe}

In \cite{BKL} it was we introduced a concept of quantum entropy and pressure which was used here.
In the classical dynamical setting the  entropy is obtained via the limit of  finite dynamical partitions.
The definition of quantum  entropy considered in \cite{BKL} is  not inspired in this procedure.  The version considered in \cite{BKL}   (and also in here)   is  more on the spirit of Rokhlin Formula (see \cite{Ma}).
This is in accordance with the point of view of considering {\it a priori} probabilities in the definition of entropy (see \cite{LMMS}).  Although at first glance it is not transparent, the classical concept of entropy is associated with an {\it a priori } measure (not an {\it a priori }  probability according to the discussion in \cite{LMMS} and \cite{BCLMS}).

In  order to make more clear the need  - in the Dynamical Systems setting - of an {\it a priori} probability consider the
 shift $\sigma$ acting on the symbolic space $K^{\mathbb{N}}$, where $K$ is a compact metric space (not necessarily countable as in \cite{BCLMS} where $K$ is the unit circle). In this case, entropy can not be defined in a natural way  by the classical procedure obtained via 
the limit of dynamical partitions. Each point $x \in K^{\mathbb{N}}$ may have an uncountable number of preimages by the shift. We consider a probability $\mu$ on $K$ which will be called the {\it a priori} probability. Entropy can be defined however, via the Ruelle operator, which in turn requires an {\it a priori} probability to be defined (in the dynamical setting this is described in \cite{LMMS}).  Entropy and the Ruelle operator are connected
in a natural and fundamental way. If $K=\{1,2,...,d\}$ is natural to consider $\mu$ as the counting measure (or the normalized counting probability) and this will give rise to the classical Ruelle operator (as described in  \cite{PP}).
The dual of the Ruelle operator plays a very important role in thermodynamic formalism (see \cite{PP}).

The time evolution considered in section \ref{gen} requires a quantum entropy with a dynamical flavor.  In the present setting, the probability $\mu$ on $M_k$  of  Section \ref{bas} will play the role of an {\it a priori probability}.  
To emphasize the differences of the various concepts of entropy we denote in this appendix the entropy of the channel $\phi_L$, when $L$ is stochastic, by $\mathfrak{h} (\phi_L)$.

Our point of view for defining the concept of  quantum entropy is in consonance with \cite{baraviera2010thermodynamic}, \cite{BLLT1}, \cite{BLLT2} and  \cite{slomczynski2003dynamical}. 
A confirmation that the entropy considered here is
in fact a concept that describes valuable information, from a dynamic point of view, is its relationship with 
Lyapunov exponents in the case of the {\it Markov model in quantum information}  (see expression \eqref{saca1}).

Section 9 in \cite{BKL} describes several distinct concepts of quantum entropy which are well known in the literature in Quantum Information (but different from ours). As mentioned there the  von Neumann entropy is not exactly of dynamical nature.

Consider  $H: M_k \to M_k$  we defined the  Ruelle operator
$\phi_H:M_k \to M_k$:
\[\rho \,\to\,\phi_H(\rho) = \int_{M_k} H(v) \rho {H(v)}^\dagger \, \dm(v),\]
 where $\mu$ is a probability on  $M_k$.
 $H$ plays the role of an Hamiltonian.

\begin{figure}[htp]
\centering
\includegraphics[width=12cm]{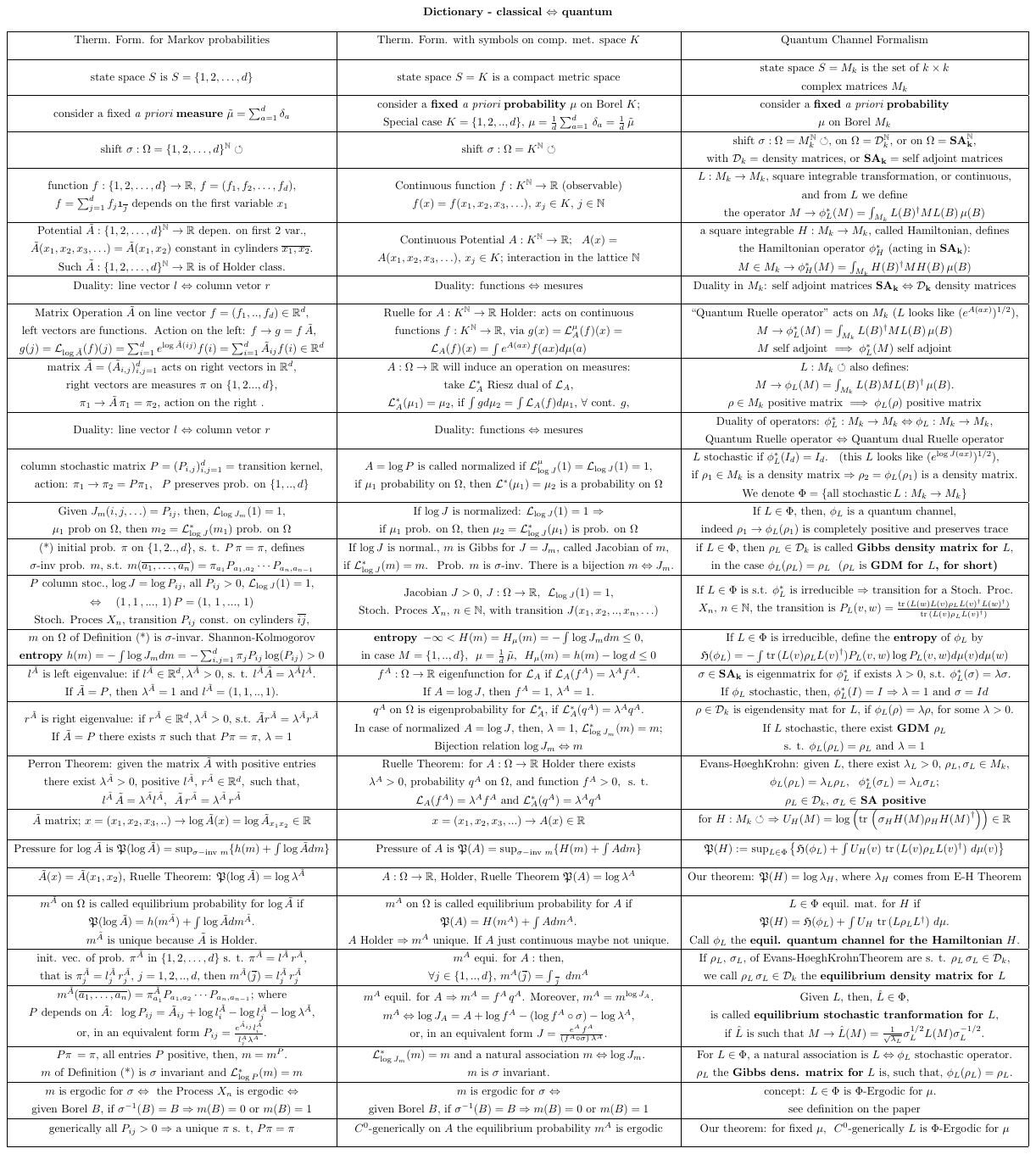}
\caption{Table 1}
\end{figure}

In \cite{BKL} it was  proved a version of the Ruelle Theorem: a variational principle  of pressure related to an eigenvalue problem for the Ruelle operator. Entropy and the Ruelle operator are also connected
in a natural and fundamental way in our  setting (see section 4 in \cite{BKL}).

In cite \cite{BKL}, given an Hamiltonian  $H: M_k \to M_k$, it is introduced the concepts of Gibbs channel and equilibrium channel for $H$. 

We present in this appendix three tables where the reader can compare different settings.
The first table, with three columns, describes the  set of correspondences (definitions and results) on the topics: in the left side of the table we consider   Markov Process (with values on a finite set), in the middle we outline results in Thermodynamic Formalism for a symbolic space with symbols in a metric space, and on the right side of the table we present the setting of Quantum Information. This table presents an outline of the sequence of results of \cite{BKL}. We believe this will help the reader to understand our reasoning.
For example: the Perron Theorem, the Ruelle Theorem, and the Evans-H{\o}egh-Krohn Theorem are on the same
line, each one in one distinct column. They play similar roles in each one of the settings.

 In  the second table we outline the main points describing  the {\it Markov model in quantum information} considered  in Section  \ref{exam}.

In the third table we make a comparison of  the results on Lyapunov exponents of the present paper with the ones in  Markov Chains and Thermodynamic Formalism.

 In the tables [BFPP] corresponds to results in \cite{benoist2017invariant}.
 
 \smallskip
 
\begin{figure}[htp]
\centering
\includegraphics[width=12cm]{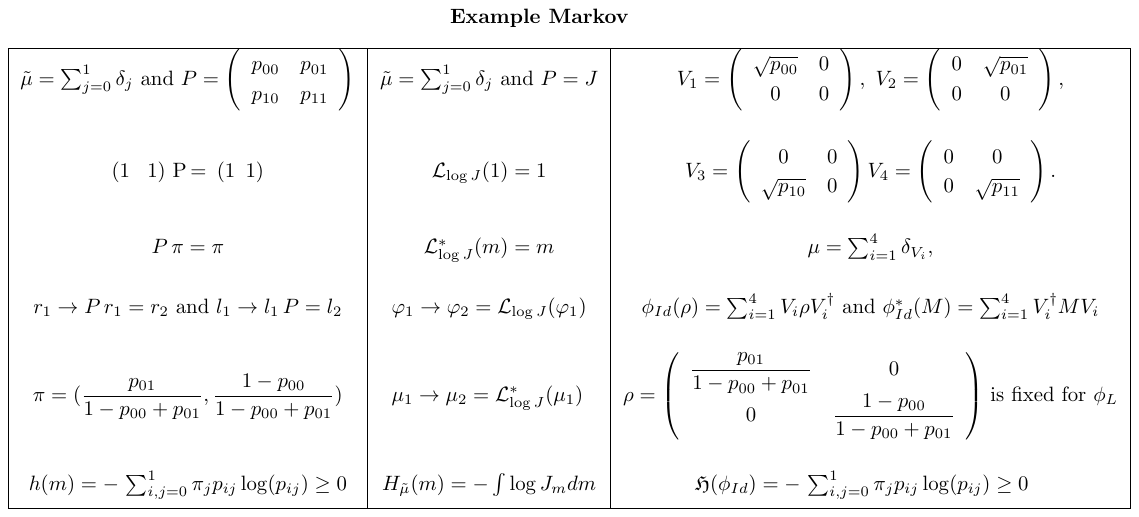}
\caption{Table 2}
\end{figure}

In order to make a comparison of our main result with the classical study of  Lyapunov exponents in Ergodic Theory, we briefly describe below some important theorems in this  long-established theory.

Given a $C^{1+\alpha}$ expanding $d$ to $1$ transformation  $T:S^1 \to S^1$, there exist a unique ergodic probability measure $ d \rho (x)$ (abs. con. with  respect to Lebesgue measure) for  $T$. The Kolmogorov entropy of this probability satisfies
\begin{equation}\label{kjl} h(\rho) = \int \log |T^\prime (x)|\, d \rho (x).
\end{equation}

This is a nice relation between entropy and Lyapunov exponent in  a smooth case (see \cite{PY}).

The Oseledet Theorem
describes the  ergodic properties of the product of random matrices (see \cite{Ma1}).

 Now, we will describe a more complex example concerning the smooth case: consider $M$ a compact manifold of dimension $2$ and a $C^\infty$ Anosov diffeomorphism $ F:M \to M$  which is volume-preserving. We denote $d F(x)$ the derivative of $F$ in $x$, which is defined for   tangent vectors $v\in T M_x$.
The tangent plane $TM_x$ is decomposed in the sum of  two bundles $E_s(x) \bigoplus E_u(x)$, where the norm of the vectors  $v$ on $E_u (x)$  contracts by the action of
$DF^{-1}$ along orbits by $F$. The bundle $E_u(x)$ is called the unstable bundle.

 We denote the Lebesgue probability by $m$, which is known to be is ergodic for such diffeomorphism $F$.
We consider the function $x \in M \to DF^{-1} (x)$ and then the random products  of the inverse of the derivative $DF^{-1}$ along orbits.

 A version of Osedelet Theorem claims in this case that there exist  $\gamma_1>0>\gamma_2$, s. t
 $$\,\lim_{n}\frac{1}{n} \log |\,[\,( D F^{-1}  (F^{n-1}(x))\,...\, D F^{-1} (x)\,]\,(v)\,|= \, \gamma_1\,\,\text{or}\,\, \gamma_2,$$
for  Lebegue  a.e.w. $x$, and $\forall v \in TM_x$. Moreover, $\gamma_1 + \gamma_2=0.$

The value $\gamma_2<0$ describes the mean value of contraction of the norm of vectors  $v \in E_u(x)$ by $DF^{-n}$, in the logarithm scale, of a random $x\in M$ under the iteration by $F^n$.

 An important result in the theory is the following:
 
 Pesin Theorem: Consider a compact  manifold $M$ of dimension $2$ and a $C^\infty$ Anosov diffeomorphism $ F:M \to M$,  which is volume preserving. Then,
\begin{equation} \label{pepe} \gamma_2 =  - h(m),
\end{equation}
 where $h(m)$ is the Kolmogorov  entropy of the Lebesgue probability $m$.

 We point out that Sections 13.1 and 13.2 in \cite{Alicki} consider results related to quantum versions of Ruelle inequality and Pesin Theorem but for a different setting.

 We will   show that our result  - for the  Lyapunov exponent of  a  Markov model in Quantum Information  - is analogous   to the Pesin Theorem. We will elaborate on this claim following the notation of \cite{benoist2017invariant}. Thus, the reader will have a synthetic view of the  chain of results in  \cite{benoist2017invariant} which were used in our work.

Take $\Omega = M_k^{\N}$, and for $n\in \N$ let $\mathcal{O}_n$ be the $\sigma$-algebra on $\Omega$ generated by the $n$-cylinder sets. An element on $\Omega$ is of the form $(\omega_1,\omega_2,...,\omega_n,...).$

We also denote
	$W_n(\omega) = L(\omega_n)...L(\omega_1).$

	Taking $O_n \in \mathcal{O}_n$ it was  defined

	$$\Pr(O_n) = \int_{O_n} \tr(W_n(\omega)\rho W_n(\omega)^*) \,\,\dm^{\otimes n}(\omega).$$

	If $\mathcal{O}$ is the smallest $\sigma$-algebra of $\Omega$ that contains all $\mathcal{O}_n$ we can extend the action of $\Pr$  to this $\sigma$-algebra.
	\smallskip

	The probability $\Pr$ on $\Omega$ defines a Stationary Stochastic Process.
\smallskip

{\bf Theorem} $(\Omega, \Pr, \theta)$ is ergodic where $\theta$ is the shift map.

 \smallskip
	The above theorem has been proved wen $L=Id$ in Lemma 4.2 in \cite{benoist2017invariant}. 
	
	\medskip

\begin{figure}[htp]
\centering
\includegraphics[width=12cm]{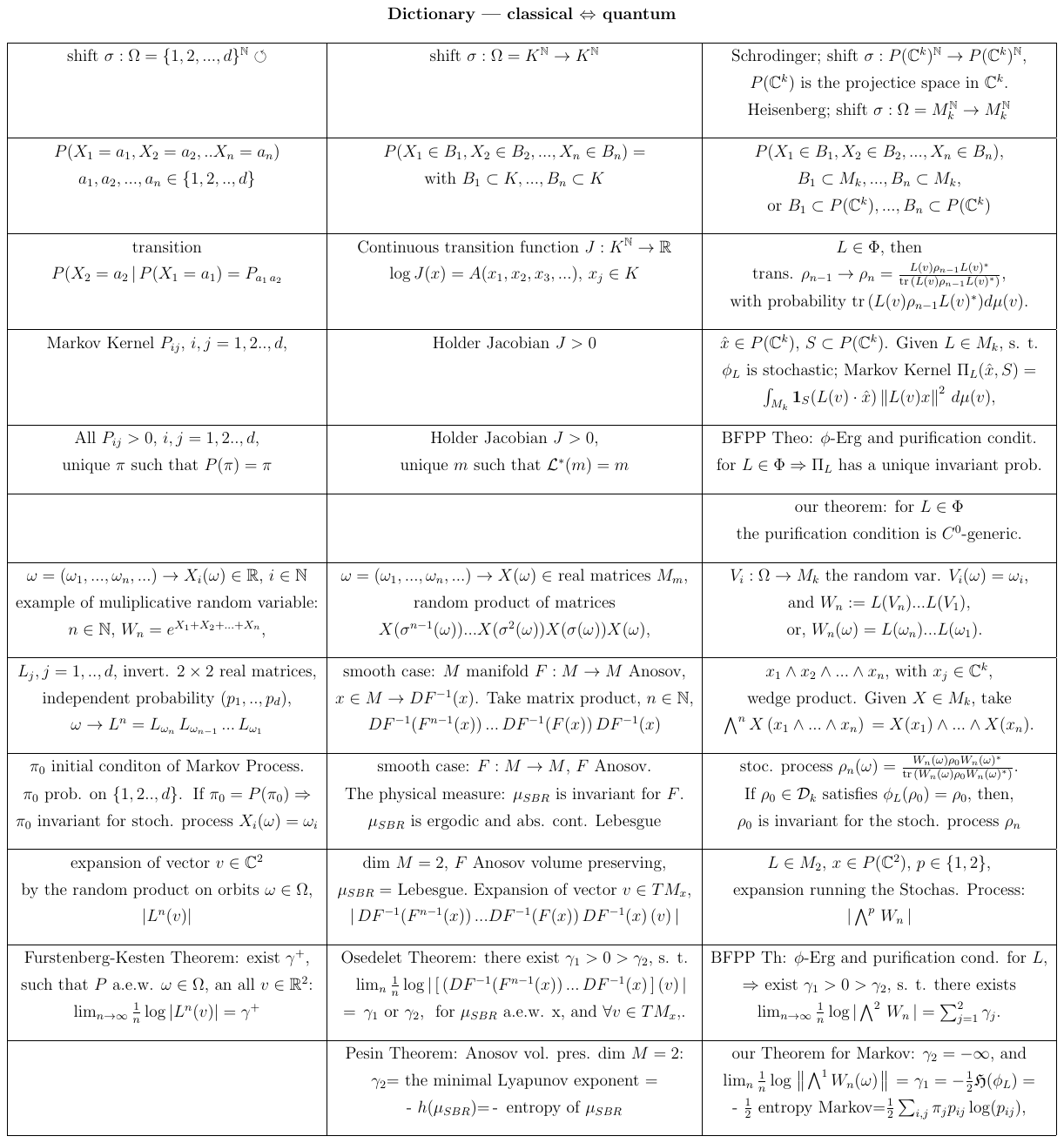}
\caption{Table 3}
\end{figure}

	Now,  we will briefly  describe the results  that are necessary for the study of Lyapunov exponents in the setting of quantum information (adapted from \cite{benoist2017invariant}).

$x_1\wedge x_2 \wedge ... \wedge x_n$, with $x_j\in \mathbb{C}^k$, denotes the classical
wedge product (an alternate form on $\mathbb{C}^k$).

One can consider an inner product
$$    \langle\, r_1\wedge r_2 \wedge ... \wedge r_n\,\,, \, \,s_1\wedge s_2 \wedge ... \wedge s_n \,\ \rangle   = \,\,\text{det}\,\, ( r_i s_j ) )_{i,j=1,2,..., n} ,
$$
and, the associated norm $ | x_1\wedge x_2 \wedge ... \wedge x_n\,|.$

Given an operator $X: \mathbb{C}^k \to \mathbb{C}^k $ we define
$ \bigwedge^n X \, : \wedge^n \mathbb{C}^k \to  \wedge^n \mathbb{C}^k $ by
$\bigwedge^n X \,( x_1\wedge x_2 \wedge ... \wedge x_n)\,= X( x_1)\wedge X(x_2) \wedge ... \wedge X(x_n).$

	\medskip

	The next theorem  is an easy adaptation of the reasoning  which was considered in \cite{benoist2017invariant} for the case $L=Id$. Below we consider the case of dimension $2$.
	
	\medskip

	{\bf Theorem.} 
		Suppose the pair $(L,\mu)$ satisfies irreducibility, the $\phi$-Erg and the purification condition. Assume also that $\int |L(v)|^2 \log |L(v)|\, \ d\mu(v)<\infty$, then, there exists numbers
		$$\infty > \gamma_1\geq \gamma_2\geq -\infty,$$
		such that, for any probability $\nu$ over $P(\mathbb{C}^k$) and any $p \in \{1,2\}$
		$$ \lim_{n \to \infty} \frac{1}{n} \log | \bigwedge^p \, W_n\,| =\sum_{j=1}^p \gamma_j,$$
		$\mathbb{P}_{\nu}$-a.s.

Consider a $2 \times 2$  stochastic matrix $P$ and the associated Markov invariant probability in $\{0,1\}^\mathbb{N}.$ 

		 In \cite{BKL} (see also \cite{baraviera2010thermodynamic}, \cite{BLLT1} and \cite{BLLT2}) it was shown that 
		 $$\mathfrak{H}(\phi_{Id})=   -\,\sum_{i,j=0}^1 \pi_j p_{ij} \log (p_{ij})=h\geq 0,$$ 
		 where $h$ is the entropy of the Markov invariant measure associated to the matrix $P$ and $\phi_{Id}$ is the channel (the Markov model in quantum information)  described in our section \ref{exam}.

		One of our main results here is for the Markov model in quantum information:
		\begin{equation} \label{olo} \gamma_1 = -\frac{1}{2}\mathfrak{H}(\phi_{Id}),
		\end{equation}
 and also that $\gamma_2= - \infty$.
 \medskip
 
 The above result is a quantum information version of \eqref{pepe}.

\end{document}